\documentclass[a4paper,11pt]{article}

\usepackage[utf8]{inputenc}
\usepackage[english]{babel}

\usepackage{amssymb,amsmath,amsthm,amsfonts}
\usepackage{mathtools}

\theoremstyle{plain}%
\newtheorem{theorem}{Theorem}[section]
\newtheorem{proposition}[theorem]{Proposition}
\newtheorem{lemma}[theorem]{Lemma}
\newtheorem{corollary}[theorem]{Corollary}

\theoremstyle{definition}
\newtheorem{definition}[theorem]{Definition}

\theoremstyle{remark}
\newtheorem{remark}[theorem]{Remark}

\usepackage{graphicx}
\usepackage{makeidx}
\usepackage{latexsym}

\newcommand{\cale}{{\mathcal E}}
\newcommand{\E}{{\mathbb E}}

\newcommand{\A}{{\mathbb A}}

\newcommand{\I}{{\mathbb I}}

\newcommand{\spec}{{\mathcal S}}

\newcommand{\RR}{{\mathbb R}}

\newcommand{\Q}{{\mathbb Q}}

\newcommand{\Pro}{{\mathbb P}}

\newcommand{\one}{{\mathbf 1}}
\newcommand{\taum}{\tau_{{\rm mix}}}

\newcommand{\n}{\|}
\newcommand{\ra}{\rightarrow}

\newcommand{\ep}{\epsilon}

\newcommand{\un}{\bigcup}

\newcommand{\iy}{\infty}

\newcommand{\bi}{\begin{itemize}}
\newcommand{\ei}{\end{itemize}}
\newcommand{\be}{\begin{enumerate}}
\newcommand{\ee}{\end{enumerate}}
\newcommand{\beq}{\begin{equation}}
\newcommand{\eeq}{\end{equation}}
\newcommand{\btm}{\begin{theorem}}
\newcommand{\etm}{\end{theorem}}
\newcommand{\bpf}{\begin{proof}}
\newcommand{\epf}{\end{proof}}
\newcommand{\bla}{\begin{lemma}}
\newcommand{\ela}{\end{lemma}}
\newcommand{\bdn}{\begin{definition}}
\newcommand{\edn}{\end{definition}}
\newcommand{\bpn}{\begin{proposition}}
\newcommand{\epn}{\end{proposition}}
\newcommand{\bcy}{\begin{corollary}}
\newcommand{\ecy}{\end{corollary}}

\DeclareMathOperator*{\Span}{Span}

\begin{document}

\title{The spectrum and convergence rates of exclusion and interchange processes on the complete graph 
}

\author{Malin P. Forsstr\"om         \and
        Johan Jonasson 
}

\date{\today}

\maketitle

\begin{abstract}
We give a short and completely elementary method to find the full spectrum of the exclusion process and a nicely limited superset of the spectrum of the interchange process (a.k.a.\ random transpositions) on the complete graph.
In the case of the exclusion process, this gives a simple closed form expression for all the eigenvalues and their multiplicities.
This result is then used to give an exact expression for the distance in \( L^2 \) from stationarity at any time and upper and lower bounds on the convergence rate for the exclusion process.
In the case of the interchange process, upper and lower bounds are similarly found.
Our results strengthen or reprove many known results about the mixing time for the two processes in a very simple way.
\end{abstract}

\section{Introduction}\label{sec:1}

Let \( G=G_n=(V,E) \) be the complete graph on $n$ vertices. The (unlabelled) {\em exclusion process} (UEP) with parameter \( \ell \) and intensity \( \alpha \) (with \( \ell \leq n/2 \) a positive integer and \( \alpha \in \RR_+ \)) on \( G \) is the continuous time Markov process \( \{X_t\}_{t \geq 0} \) on the set \( {V \choose \ell} \) of \( \ell \)-element subsets of \( V \), defined by taking its generator \( \Q = \Q^{(n,\ell,\alpha)} = [q_{J,J'}]_{J,J' \in {V \choose \ell}} \) as
\[q_{J,J'} = \left\{ \begin{array}{ll} -\alpha \ell(n-\ell),& J=J' \\ \alpha, & |J \Delta J'|=2 \\ 0, & \mbox{otherwise} \end{array} \right.\]
Clearly \( \Q \) is symmetric.
We usually think of this process as having either a black or a white ball at each vertex \( v \in V \), letting the state denote the set of \( \ell \) vertices where there is a black ball.
For each edge \( e = \{u,v\} \in E \), we can associate a Poisson clock of intensity \( \alpha \) such that whenever the clock rings, the two balls at \( u \) and \( v \) switch positions.
(Since the black balls are not distinct, this means that the process jumps to a new state only when the Poisson clock of an edge with one black and one white ball rings.)
We take \( X_0 \) to be an arbitrary but fixed state \( J \).

The {\em labelled exclusion process} (LEP) with parameters \( \ell \) and \( \alpha \) is the same process with the exception that we replace the black balls with \( \ell \) {\em distinct}
balls, with labels (or colors (not white) if you like) \( 1,\ldots,\ell \). Here we may also take \( n/2 < \ell \leq n \).
The state space is now the set \( (V)_\ell \) of ordered \( \ell \)-tuples of distinct elements of \( V \). For \( x=(v_1,\ldots,v_\ell) \in (V)_\ell \), we will think of \( v_j \) as the position
of ball \( j \).
Obviously \( |(V)_\ell| = (n)_\ell = n(n-1)\ldots (n-\ell+1) \) and the generator \( \Q = \Q^{(n,\ell,\alpha)} = [q_{x,x'}]_{x,x' \in (V)_\ell} \) is given by
\( q_{x,x} = -\alpha(\ell(n-\ell)+{\ell \choose 2}) \), \( q_{x,x'} = \alpha \) whenever \( x \) and \( x' \) differ for exactly one ball or \( x' \) can be obtained from \( x \) by interchanging two of its elements, and \( q_{x,x'}=0 \) otherwise.
Again \( \Q \) is symmetric.
The special case \( \ell = n \) makes \( (V)_\ell \) the set of permutations of \( n \) balls, in which case the process is also known as the {\em interchange process} or
{\em random transpositions} on \( G \).

The alert reader will have spotted an ambiguity in our notation: we have used \( \Q \) for the generator of two different processes.
However, this should be no problem, since it will always be clear which one is under consideration.

The spectrum of the UEP is known, see e.g.~\cite{BCN,DS,M}:
\btm \label{ta}
Let \( \Q^{(n,\ell,\alpha)} \) be the generator of the UEP with parameters \( \ell \) and \( \alpha \).
Then the eigenvalues of \( -\Q^{(n,\ell,\alpha)} \) are
\[0,\alpha n,2\alpha(n-1),3\alpha(n-2),\ldots,\ell\alpha(n-\ell+1)\]
with multiplicities
\[1,n-1,{n \choose 2}-n,{n \choose 3}-{n \choose 2},\ldots,{n \choose \ell}-{n \choose \ell-1}\]
respectively.
\etm

To the best of our knowledge, most previous proofs of Theorem~\ref{ta} require a fair deal of background knowledge, whereas our short proof requires nothing beyond
standard undergraduate linear algebra.
Moreover, our method generalizes fairly easily to the LEP to find a nicely limited set which contains the full spectrum of that process. This spectrum can be understood using representation theory (see e.g. Wimmer~\cite{W}), but this requires much more theory.
%(Of course, the second eigenvalue is well known to be \( \alpha n \) for both processes, see e.g.\ Diaconis and Shahshahani \cite{DS}. Indeed, this follows from a famous conjecture of Aldous that the spectral gap of the interchange process on any graph \( G \) is the same as for random walk on \( G \), a result now known to be true, see Caputo, Liggett and Richthammer \cite{CLR}.)

To present our main result for the spectrum of the LEP, we let \( \A_n \), \( n=1,2,\ldots \) be the adjacency matrix for the Cayley graph of the symmetric group on \( n \) items generated by the transpositions,
i.e.\ the graph which for vertex set has all the \( n! \) permutations and an edge between \( u \) and \( v \) iff \( u \) and \( v \) differ by exactly one transposition.
The eigenvalues \( \mu_j \), \( 1 \leq j \leq n! \), of \( \A_n \) relate to the eigenvalues \( \lambda_j \)  of \( -\Q^{(n,n,\alpha)} \), (i.e.\ the generator for the LEP with \( \ell=n \)) by
\[\lambda_j=\alpha\left({n \choose 2}-\mu_j\right).\]
For any square matrix \( {\mathbb B} \), write \( \spec({\mathbb B}) \) for the set of eigenvalues of \( {\mathbb B} \).
Let \( \cale^n_0 = \{0\} \), \( \cale^n_1=\{0,\alpha n\} \) and inductively for \( k=1,2,\ldots,n -2\),
\[\cale^n_{k+1} = \cale^n_k \un \left( \alpha\left(n(k+1)-{k+1 \choose 2} - \spec(\A_{k+1})\right) \right) \]
and \( \cale^n_n = \cale^n_{n-1} \). Here, for \(x,y \in \mathbb{R}\) and \( A \subset \mathbb{R} \), we use \( x (y - A) \) to denote the set of all real numbers which can be written as \( x(y-a) \) for some \( a \in A \).

\btm \label{tb}
Let \( \Q^{(n,\ell,\alpha)} \) be the generator of the LEP with parameters \( n \), \( \alpha \) and \( \ell \). Then with \( n \) fixed and \( \cale^n_0,\ldots,\cale^n_n \) as defined above,
\( \spec(-Q^{(n,\ell,\alpha)}) \)
is increasing in \( \ell \) and
\[\spec(-\Q^{(n,\ell,\alpha)}) \subset \cale^n_\ell.\]
Moreover, \( \spec(-\Q^{(n,n,\alpha)}) \) is contained in the largest subset of \( \cale^n_n \) that is symmetric around \( \alpha{n \choose 2} \) or, equivalently
\( \spec(\A_n) \) is symmetric and is  contained in the largest subset of \( {n \choose 2} - \cale^n_n \) that is symmetric around \( 0 \).

Furthermore, if \( 1 \leq j \leq \min{ (\ell,n-\ell)} \), then the total multiplicity of the eigenvalues \( \lambda \) such that \( \alpha j(n-j+1) \leq \lambda < \alpha (j+1)(n-j) \), does not exceed \( (n)_j {\ell \choose j} \).
Also, for all \( \ell \), the multiplicity of the eigenvalue \( \alpha n \) is exactly \( \ell(n-1) \).
\etm

Theorem~\ref{tb} can be used recursively on \( n \) and \( \ell \) to find supersets of \( \spec(-\Q^{(n,\ell,\alpha)}) \): having found supersets of \( \cale^j_k \) for \( j<n \) and \( k \leq j \), we find supersets of
\( \spec(\A_j) \) for all \( j<n \) and then the \( \cale^n_k \):s.

\begin{remark}
 Note that it is obvious that \( \spec(-\Q^{(n,n,\alpha)}) = \spec(-\Q^{(n,n-1,\alpha)}) \).
Note also that the eigenvalues for \( \A_k \) are symmetrically spread out between \( -\alpha{k \choose 2} \) and \( \alpha{k \choose 2} \). As a consequence, for \( \ell=o(\sqrt{n}) \), the sets \( \cale^n_0,\ldots,\cale^n_\ell \) are disjoint and for \( \ell = o(n) \), the spread-outs of \( \cale^n_0,\ldots,\cale^n_\ell \) are of smaller order than their centers.
\end{remark}

The results of Theorems~\ref{ta} and~\ref{tb} have profound consequences for the time taken for these processes to come close to uniformity.
Common ways to quantify the distance between two probability measure are by the \( L^p \)-norm or the total variation norm.
Let \( \pi \) be a probability measure on a finite space \( S \). If \( \nu \) is a signed measure on \( S \), then we define the \( L^p(\pi) \) norm of \( \nu \) for \( p \geq 1 \) by
\[\n \nu \n_p^p = \E_\pi\left[\left|\frac{\nu(X)}{\pi(X)}\right|^p\right] = \sum_{s \in S}\left|\frac{\nu(s)}{\pi(s)}\right|^p\pi(s).\]
For a probability measure \( \mu \) on \( S \), the \( L^p \)-distance from \( \mu \) to \( \pi \) is then defined as \( \n \mu-\pi \n_p \).
By H\"older's inequality, \( \n \mu-\pi \n_p \) is increasing in \( p \). The total variation distance is defined as
\[\n \mu - \pi \n_{TV} = \frac12 \n \mu - \pi \n_1.\]

To define what we mean by the mixing time for a Markov chain, let \( \{X_t\} \) be a Markov chain on \( S \) having stationary distribution \( \pi \) and let \( \Pro_{x_0} \) be the underlying probability measure when starting from \( X_0=x_0 \).
Then the mixing time of \( \{X_t\} \) is defined for any \( \ep \in (0,1) \) as
\[\taum(\ep) = \inf\{t \colon \max_{x_0}\n\Pro_{x_0}(X_t \in \cdot)-\pi\n_{TV} \leq \ep\}.\]
For \( p > 1 \), the \( L^p \)-mixing time is defined as
\[\tau_p(\ep)=\inf\{t \colon \max_{x_0}\n\Pro_{x_0}(X_t \in \cdot) - \pi \n_p \leq 2\ep\}.\]
Hence \( \taum(\ep) = \tau_1(\ep)  \leq \tau_p(\ep) \) and \( \tau_p(\ep) \) is increasing in \( p \).
One standard is to work with \( p=2 \), which is the norm that is most naturally associated to the spectrum of the Markov chain.
Note that
\[\n\mu-\pi\n_2^2 = \sum_{s \in S} \frac{(\mu(s)-\pi(s))^2}{\pi(s)}\]
which in case \( \pi \) is uniform becomes
\[\n\mu-\pi\n_2^2 = |S|\sum_{s \in S}(\mu(s)-\pi(s))^2.\]

Often results on mixing times are very precise in an asymptotic sense as the size of the state space goes to infinity.
In such cases, we are in fact considering a sequence of Markov chains \( \{X^n_t\} \) on state spaces \( S^n \) such that \( |S^n| \ra \iy \)
and we try to express \( \taum^n(\ep) \) or \( \tau_p^n(\ep) \) in terms of \( n \).
Usually there is an obvious natural way to define the \( \{X^n_t\} \):s and the \( S^n \):s.
In our case we will simply let the number of vertices, \( n \), grow.

Our main results are the following. We set \( \alpha \) to \( 2/n^2 \) in order to get the standard case of one state change per time unit, but the results easily generalize to arbitrary
\( \alpha \) if you like. By symmetry, \( \n \Pro_{x_0}(X_t \in \cdot) - \pi \n_p \) does not depend on \( x_0 \), neither for the UEP nor the LEP, so \( x_0 \) has been dropped from the notation.

\btm \label{tc}
Let \( \{X_t\}_{t \geq 0} \) be the unlabelled exclusion process with \( n \) balls in total and \( \ell=\ell(n) \) black balls and set \( \alpha=2/n^2 \).
Then for any \( \ell \leq n/2 \),
\[\n\Pro(X_t \in \cdot)-\pi\n_2^2 = \sum_{i=1}^{\ell} \left( {n \choose i} - {n \choose i-1}\right)e^{-4i(n-i+1)t/n^2}.\]
As a consequence, writing \( t=(1/4)n\log(n-1)+cn \) for a constant \( c \),
\[e^{-2c} \leq \n\Pro(X_t \in \cdot) - \pi\n_2  \leq 2e^{-2c},\]
where the upper bound holds for \( c \geq 0 \) and sufficiently large \( n \).
In particular for all \( \ell \) and all \( \ep \in (0,1) \),
\[\tau_2(\ep) = \frac14 n\log n + C(\ep)n\]
for a a constant \( C(\ep) \) depending on \( \ep \).
\etm

\begin{remark}
 Lacoin and Leblond \cite{LL} proved that
\[
 \taum(\ep) = \left(1+o(1)\right)\frac12 n\log \min(\ell,\sqrt{n}).
\]
Our result confirms the upper bound of this result for \( \ell \geq \sqrt{n} \).
For \( \ell < \sqrt{n} \) we note that there is a significant difference between \( \taum(\ep) \) and \( \tau_2(\ep) \) which comes from the fact that the \( L^2 \) norm is much less forgiving about
any remaining traces of the starting state.
However, to establish the upper bound on \( \taum \) for \( \ell<\sqrt{n} \) can be readily done by a straightforward coupling argument.

For $\ell \geq \sqrt{n}$, Theorem \ref{tc} shows that $\taum(\ep) \leq (1+o(1))(1/4)n\log n$. For a matching lower bound, consider the number of black balls that at time $t$ are in
positions that had a black ball at time $0$.
Taken together, these facts establish that there is a cutoff in total variation at time $(1/4)n\log n$.
\end{remark}

\begin{remark}
One may analyze the exact expression for $\n \Pro(X_t \in \cdot)-\pi\n_{2}^2$ in Theorem \ref{tc} asymptotically as $n \ra \iy$.
Using essentially the same computations as below in the proof of Theorem \ref{tc}, one fairly easily finds that if also $\ell \ra \iy$, then
$\n \Pro(X_t \in \cdot)-\pi\n_{2}^2 = (1+o(1))(e^{e^{-4c}}-1)$.
In case $\ell$ stays constant, the asymptotic expression on the right hand side becomes $(1 + o(1))\sum_{i=1}^{\ell} e^{-4ci}/i!$, or equivalently  $(1+o(1))(e^{e^{-4c}}-1-\sum_{i=\ell+1}^\iy e^{-4ci}/i!)$.
\end{remark}

\btm \label{td}
Let \( \{X_t\}_{t \geq 0} \) be the labelled exclusion process with \( n \) balls in total, \( \ell = \ell(n) \leq (1-\varepsilon)n/2 \) labelled balls for some \( \varepsilon > 0 \) and  \( \alpha = 2/n^2 \).
Then for \({ t=(1/4)n\log(\ell(n-1)) + cn} \),
\[e^{-2c} \leq \n \Pro(X_n \in \cdot)-\pi \n_2 \leq 2e^{-2c},\]
where the upper bound applies for \( c \geq 0 \) and sufficiently large \( n \).
In particular for all \( \delta \in (0,1) \),
\[\tau_2(\delta) = \frac14 n \log(\ell n) + C(\delta)n.\]
\etm

\begin{remark}
As for the UEP, straightforward probabilistic arguments prove that for \( \ell = o(n) \), \( \taum(\ep) = (1/2)(1+o(1))\log \ell \), so here we have a significant difference
between \( \taum(\ep) \) and \( \tau_2(\ep) \) for all \( \ell = o(n) \).
Theorem~\ref{td} does not establish the well-known fact that for \( \ell=n \), \( \tau_2(\ep) = (1+o(1))(1/2)n\log n \), see \cite{DS}. However, it at least follows from Theorem~\ref{tb} that
the relaxation time is \( n/2 \).
The lower bound for \( \ell=n \) follows easily by probabilistic arguments, even for total variation, simply considering the number of labelled balls that are still in their
starting positions, see \cite{DS}.
\end{remark}

%\begin{remark}
%Our proof of Theorem~\ref{td} can easily be extended to any \( \ell \leq \beta n \)  for any \( \beta \in ( 0 ,1/2 ) \) by simply adjusting the constants used.
%\end{remark}

An outline of the remainder of the paper is as follows. Theorems~\ref{ta} and~\ref{tb} are proved in Section 2 and Section 3 respectively and the short proofs of
Theorems~\ref{tc} and~\ref{td} are then given in Section 4.

\section{Spectrum for the unlabelled exclusion process}\label{sec:2}

In order to lighten the notation, we will assume here that \( \alpha=1 \); the generalization to the case of arbitrary \( \alpha \) is trivial.
We will also consider \( n \) fixed. Hence we will simply write \( \Q^{(\ell)} \) for \( \Q^{(n,\ell,1)} \).
We will start by showing that if \( \phi \colon {V \choose \ell} \ra \RR \) is an eigenvector for \( \Q^{(\ell)} \), then it can be {\em lifted} to an eigenvector for \( \Q^{(\ell+1)} \).
For a set \( J \in {V \choose \ell+1} \) and \( v \in V \), write \( J_v \coloneqq J \setminus \{v\} \).
For an edge \( e=\{u,v\} \), let \( J_e \) be the state one gets from flipping the balls at \( u \) and \( v \).

\bdn \label{da}
For a function \( f\colon{V \choose \ell} \ra \RR \), let the lift of \( f \) on \( {V \choose \ell+1} \) be given by
\[\hat{f}(J) = \sum_{j \in J}f(J_j).\]
\edn

\bla \label{la}
Assume that \( \phi \) is a nonzero eigenvector of \( -\Q^{(\ell)} \) with corresponding  eigenvalue \( {\lambda \leq \ell(n-\ell+1)}\).
Then \( \hat{\phi} \) is a nonzero eigenvector of \( -\Q^{(\ell+1)} \) for the same eigenvalue.
\ela

\bpf
The crucial observation is that for any \( J \in {V \choose \ell+1} \), we have
\[\sum_{e \in E} \hat{\phi}(J_e) = \sum_{e \in E} \sum_{j \in J_e} \phi((J_e)_j) = \sum_{e \in E} \sum_{j \in J} \phi((J_j)_e).\]
This implies that
\begin{align*}
-\Q^{(\ell+1)}\hat{\phi}(J)
&= \sum_{e \in E}(\hat{\phi}(J)-\hat{\phi}(J_e))
= \sum_{e \in E} \sum_{j \in J} (\phi(J_j)
- \sum_{e \in E} \sum_{j \in J_e} \phi((J_e)_j)
\\&= \sum_{e \in E} \sum_{j \in J} (\phi(J_j)-\phi((J_j)_e)
= \sum_{j \in J} -\Q^{(\ell)}\phi(J_j)
= \sum_{j \in J} \lambda \phi(J_j) \\
&= \lambda \hat{\phi}(J).
\end{align*}

This proves that \( \hat{\phi} \) is either an eigenvector of the desired form, or the zero vector.
To rule out the second possibility, we observe that if this were the case, then by definition, for any \( J \in {V \choose \ell+1} \),
\[\sum_{j \in J} \phi(J_j) = 0.\]
A given \( K \in {V \choose \ell} \) gives rise to a term in the left hand side sum for all \( J \in {V \choose \ell+1} \) such that \( J \supset K \). It follows that for any such \( J\), we have
\[\phi(K)+\sum_{K'\subset J \colon |K \Delta K'|=2}\phi(K') = 0.\]
Summing over \( J \supset K \), we get
\[(n-\ell)\phi(K) + \sum_{K' \in {V \choose \ell}\colon |K \Delta K'|=2}\phi(K') = 0.\]
We recognize the sum above as \( (\ell(n-\ell)\I_\ell+\Q^{(\ell)})\phi(K) \), where \( \I_\ell \) is the identity matrix of dimension \( \binom{n}{\ell} \). Since this holds for all \( K \), it follows that the system of equations
\[((\ell+1)(n-\ell)\I_{\ell}+\Q^{(\ell)})\phi = 0\]
holds.
Since by assymption, the eigenvalue \( \lambda \) corresponding to \( \phi \) satisfies \( -\lambda \geq -  \ell(n-\ell+1) \), no eigenvalue of \( \Q^{(\ell)} \) is smaller than \( -\ell(n-\ell+1) \) and \( \ell<n/2 \), \( 0 \) is not an eigenvalue of \( (\ell+1)(n-\ell)\I_\ell+\Q^{(\ell)} \) and hence \( \phi \equiv 0 \) is the unique solution.
This contradicts that \( \phi \) is a nonzero eigenvector.
\epf

Define the usual inner product on \( L^2\bigl({V \choose \ell},\pi\bigr) \) by
\[\langle \phi,\psi \rangle = \E_\pi[\phi(X)\psi(X)] = {n \choose \ell}^{-1}\sum_{J \in {V \choose \ell}}\phi(J)\psi(J)\]
and say that \( \phi \) and \( \psi \) are orthogonal if their inner product is \( 0 \).

\bla \label{lb}
Assume that \( 1 \leq k \leq (n/2)-1 \) and that \( \phi,\psi\colon{V \choose k} \ra \RR \) are orthogonal eigenvectors of \( \Q^{(k)} \).
Then \( \hat{\phi} \) and \( \hat{\psi} \) are orthogonal eigenvectors of \( \Q^{(k+1)} \).
\ela

\bpf
Assume that \( \phi,\psi\colon{V \choose k} \ra \RR \) are orthogonal eigenvectors of \( \Q^{(k)} \). By Lemma~\ref{la}, \( \hat{\phi} \) and \( \hat{\psi} \) are eigenvectors of  \( \Q^{(k+1)} \). As eigenvectors of a symmetric matrix, they must be orthogonal unless they are eigenvectors for the same eigenvalue
\( \lambda \), so assume that this is the case.
Then
\begin{align*}
{n \choose k+1} \langle \hat{\phi},\hat{\psi} \rangle &= \sum_{J \in {V \choose k+1}} \sum_{K \in {V \choose k}\colon K \subset J} \sum_{K' \in {V \choose k}\colon K' \subset J} \phi(K)\psi(K') \\
&= \sum_{K \in {V \choose k}}\left( (n-k)\phi(K)\psi(K) + \sum_{K'\colon |K \Delta K'|=2} \phi(K)\psi(K') \right) \\
&= \sum_{K \in {V \choose k}} \phi(K) \left( (k+1)(n-k)\psi(K) - k(n-k)\psi(K) + \sum_{K'\colon |K \Delta K'|=2}\psi(K') \right) \\
&= \sum_{K \in {V \choose k}} \phi(K) \left( (k+1)(n-k)\psi(K) + \Q^{(k)}\psi(K) \right) \\
&= \sum_{K \in {V \choose k}} \phi(K) \left( (k+1)(n-k)\psi(K) - \lambda\psi(K) \right) \\
&= ((k+1)(n-k)-\lambda) \sum_{K \in {V \choose k}} \phi(K)\psi(K) \\
&= 0
\end{align*}
where we everywhere sum over \( K' \in {V \choose k} \) and the final equality uses that \( \phi \) and \( \psi \) are orthogonal.

\epf

With these results at hand, we are ready to prove Theorem~\ref{ta}.

\begin{proof}[Proof of Theorem~\ref{ta}]
This will be done with induction over \( \ell \).
The result is trivial for \( \ell=0 \) and well-known for \( \ell=1 \).
Assume now that it also holds for \( \ell=2,3,\ldots,k \), in particular that \( -\Q^{(k)} \) has the eigenvalues \( 0,n,2(n-1),\ldots,k(n-k+1) \)
of respective multiplicities \( 1,n-1,{n \choose 2}-n,\ldots,{n \choose k}-{n \choose k-1} \).
Since \( \Q^{(k)} \) is symmetric, we can find an orthogonal set of eigenvectors \( \phi_1,\ldots,\phi_{{n \choose k}} \).
By Lemmas~\ref{la} and~\ref{lb}, \( \hat{\phi_1},\ldots,\hat{\phi}_{{n \choose k}} \) is an orthogonal set of eigenvectors
of \( \Q^{(k+1)} \) for the same eigenvalues.

Let \( H \coloneqq Span (\hat{\phi_1},\ldots,\hat{\phi}_{{n \choose k}}) \). It then only remains to prove that any vector in the orthogonal complement \( H^\perp \) of \( H \) is an eigenvector of  \( - \Q^{(k+1)} \) with eigenvalue \( (k+1)(n-k) \).
%
%
%we can find further vectors \( \psi_{{n \choose k}+1},\ldots,\psi_{{n \choose k+1}} \) that are pairwise orthogonal and all eigenvectors for the eigenvalue \( (k+1)(n-k) \).
%
To see this, note first that for any  vector \( \psi \in H^\perp \) and any  \( I \in {V \choose k} \), we must have
\begin{equation} \label{eb}
\sum_{i \not \in I} \psi(I \cup \{i\}) = 0.
\end{equation}
Also, if we spell out the equation \( -\Q^{(k+1)}\psi = \lambda \psi \) for some  \( J \in {V \choose k+1} \), we get
\begin{equation} \label{ea}
(k+1)(n-k-1)\psi(J) - \sum_{K \in {V \choose k+1} \colon |K \Delta J|=2} \psi(K) = \lambda \psi(J).
\end{equation}
If we use~\eqref{eb}, we sum in the previous equation becomes \( -(k+1)\psi(J) \) and the system of equations simply becomes that for each \( J \),
\begin{equation} \label{ec}
(k+1)(n-k)\psi(J)=\lambda\psi(J).
\end{equation}
Obviously this cannot hold for a nonzero \( \psi \) unless \( \lambda = (k+1)(n-k) \) and provided that this is so, then {\em any} nonzero vector \( \psi \in \RR^{{n \choose k+1}} \) satisfies~\eqref{ec} for all \( J \).
Since~\eqref{eb} imposes \( {n \choose k} \) linear restrictions, it follows that when  \( \lambda=(k+1)(n-k) \),  we can find \( {n \choose k+1}-{n \choose k} \) pairwise orthogonal vectors \( \psi \in \mathbb{R}^{\binom{n}{k+1}} \) solving~\eqref{eb}, and hence also~\eqref{ea}.
\end{proof}

\section{Spectrum for the labelled exclusion process}\label{sec:3}

As in the previous section, our notation of the generators will be \( \Q^{(\ell)} \) and we assume that \( \alpha=1 \).
(Note that then the \( \cale^n_k \):s contain only integer values.)
In analogy with the UEP, we will need to lift a function \( f\colon (V)_k \ra \RR \) to a function on \( (V)_{k+1} \).
However since we can now identify the balls, lifts are in fact more straightforward; define for each \( i=1,2,\ldots,k+1 \), \( f^i\colon (V)_{k+1} \ra \RR \) as
\[f^i(v_1,\ldots,v_{i-1},v_i,v_{i+1},\ldots,v_{k+1}) = f(v_1,\ldots,v_{i-1},v_{i+1},\ldots,v_{k+1}),\]
In other words, \( f^i \) is derived from \( f \) by simply ignoring the position of the \( i \)'th labelled ball.
It is then obvious that if \( \phi\colon (V)_k \ra \RR \) is an eigenvector of \( \Q^{(k)} \) for the eigenvalue \( \lambda \), then \( \phi^i \) is an eigenvector of \( \Q^{(k+1)} \)
for the same eigenvalue and that if \( \phi,\psi\colon (V)_k \ra \RR \) are orthogonal, then so are \( \phi^i \) and \( \psi^i \).
(Here of course the inner product is defined in complete analogy with the UEP.)

\begin{proof}[Proof of Theorem~\ref{tb}]
Since for \( \ell \leq 1 \), there is no difference between the UEP and the LEP, we know that Theorem~\ref{tb} holds for all \( n \) and \(\ell\leq 1 \).
Assume for induction that for some fixed \( n \), the result  holds for  \( \ell=0,1,\ldots,k \). It then suffices to prove the result for  \( \ell=k+1 \).
Note that the induction hypothesis tells us that \( \spec(-\Q^{(k)}) \subseteq \cale^n_k \) and hence in particular contains only integers and ranges at most from \( 0 \) to \( nk \).

Let \( 0=\lambda_1<\lambda_2 \leq \ldots \leq \lambda_{(n)_k} \) be the eigenvalues of \( \Q^{(k)} \) and \( \phi_1,\ldots,\phi_{(n)_k} \) be a corresponding orthogonal set of eigenvectors.
Then for any \( {i \in \{1,\ldots,k+1\}} \), \( \phi^i_1,\ldots,\phi^i_{(n)_k} \) are orthogonal eigenvectors of \( \Q^{(k+1)} \) for the same eigenvalues.
In particular any eigenvalue of \( \Q^{(k)} \) is also an eigenvalue for \( \Q^{(k+1)} \).

We claim that a vector \( \phi \) is orthogonal to all vectors in the span of \( \{\phi^i_j\}_{1 \leq i \leq k+1,1 \leq j \leq (n)_k}\) if and only if for all \( 1 \leq i \leq k+1 \) and all \( (v_1,\ldots,v_{i-1},v_{i+1},\ldots,v_{k+1}) \in (V)_k \),
\begin{equation} \label{ee}
\sum_{v\in V \colon v \not \in \{v_1,\ldots,v_{i-1},v_{i+1},\ldots,v_{k+1}\}}\phi(v_1,\ldots,v_{i-1},v,v_{i+1},\ldots,v_{k+1}) = 0.
\end{equation}

To see that the if direction holds, note simply that if \( \phi \) satisfies~\eqref{ee}
for all \( 1 \leq i \leq k+1 \) and all \( (v_1, \ldots, v_{i-1},v_{i+1}, \ldots, v_{k+1} ) \), then for all \( \phi_j^i \in \{ \phi_j^i\}_{1 \leq i \leq k+1,\, 1 \leq j \leq (n)_k}\) we have
\begin{equation}\label{mb}
\begin{split}
(n)_k \langle \phi_j^i, \phi \rangle
&=
\sum_{\mathrlap{(v_1, \ldots, v_{k+1}) \in V_{k+1}}}
\phi(v_1, \ldots, v_{k+1}) \phi_j^i (v_1, \ldots, v_{k+1} )
\\&=
\sum_{\mathrlap{(v_1, \ldots, v_{i-1}, v_{i+1}, \ldots, v_{k+1}) \in V_{k} \atop v \in V \backslash \{ v_1, \ldots, v_{i-1}, v_{i+1}, \ldots, v_{k+1}\}}} \phi_j^i(v_1, \ldots, v_{i-1}, v, v_{i+1}, \ldots, v_{k+1}) \phi (v_1, \ldots, v_{i-1}, v,  v_{i+1}, \ldots, v_{k+1})
\\&=
\sum_{\mathrlap{(v_1, \ldots, v_{i-1}, v_{i+1}, \ldots, v_{k+1}) \in V_{k} \atop v \in V \backslash \{ v_1, \ldots, v_{i-1}, v_{i+1}, \ldots, v_{k+1}\}}}
\phi_j (v_1, \ldots, v_{i-1}, v_{i+1}, \ldots, v_{k+1})
\phi(v_1, \ldots, v_{i-1}, v, v_{i+1}, \ldots, v_{k+1}) .
\\&
=
 \sum_{\mathrlap{{(v_1, \ldots, v_{i-1}, v_{i+1}, \ldots, v_{k+1}) \in V_{k} }}}
 \phi_ j(v_1, \ldots, v_{i-1},   v_{i+1}, \ldots, v_{k+1})
\sum_{\mathrlap{  v \in V \backslash \{ v_1, \ldots, v_{i-1}, v_{i+1}, \ldots, v_{k+1}\}}}
\phi  (v_1, \ldots, v_{i-1}, v,  v_{i+1}, \ldots, v_{k+1}).
\end{split}
\end{equation}
As the second sum is zero by assumption, it follows that \( \langle \phi, \phi_j^i \rangle = 0 \). As this holds for all \( \phi_j^i \in \{ \phi_j^i\}_{1 \leq i \leq k+1,\, 1 \leq j \leq (n)_k}\), then clearly \( \phi \) is orthogonal with all \( \psi \in   \Span \{ \phi_j^i\}_{1 \leq i \leq k+1,\, 1 \leq j \leq (n)_k} \).

For the other direction, suppose that \( \phi \) is orthogonal to all \(  \{ \phi_j^i\}_{1 \leq i \leq k+1,\, 1 \leq j \leq (n)_k} \), i.e. that for any \( 1 \leq i \leq k+1 \) and \( 1 \leq j \leq(n)_k \), \( \langle \phi_j^i,\phi \rangle= 0\). For any \( (v_1, \ldots, v_{i-1}, v_{i+1}, \ldots, v_{k+1} ) \in V_k \), define
\[
 \psi_i (v_1, \ldots, v_{i-1}, v_{i+1}, \ldots, v_{k+1} ) \coloneqq
\sum_{\mathrlap{  v \in V \backslash \{ v_1, \ldots, v_{i-1}, v_{i+1}, \ldots, v_{k+1}\}}}
\phi  (v_1, \ldots, v_{i-1}, v,  v_{i+1}, \ldots, v_{k+1}).
\]
Then clearly \( \psi_i \colon V_k \to \mathbb{R} \). As \( \{ \phi_j \} \) spans the set of all real valued functions from \( V_k \) to \( \mathbb{R} \), it follows that \( \psi_i \in \Span \{ \phi_j\}_{ 1 \leq j \leq (n)_k} \). Using~\eqref{mb}, it follows that
\begin{align*}
0 &= (n)_{k+1} \langle \phi_j^i, \phi \rangle
\\&=  \sum_{\mathrlap{(v_1, \ldots, v_{i-1}, v_{i+1}, \ldots, v_{k+1}) \in V_{k} }}
 \phi_ j(v_1, \ldots, v_{i-1},   v_{i+1}, \ldots, v_{k+1}) \quad
 \psi_i (v_1, \ldots, v_{i-1}, v_{i+1}, \ldots, v_{k+1})
\\&=  (n)_k \langle  \phi_ j,  \psi_i \rangle.
\end{align*}
 As  \( \psi_i \in   \Span \{ \phi_j\}_{1 \leq j \leq (n)_k} \) and this holds for all \( j \), we must have that \( \psi_i \equiv 0 \), or equivalently,
\[
\sum_{\mathrlap{\mathclap{  v \in V \backslash \{ v_1, \ldots, v_{i-1}, v_{i+1}, \ldots, v_{k+1}\}}}}
\phi  (v_1, \ldots, v_{i-1}, v,  v_{i+1}, \ldots, v_{k+1}) = 0
\]
for all \( (v_1, \ldots, v_{i-1}, v_{i+1}, \ldots, v_{k+1} ) \in V_k \).

Assume now that \( \phi \) is an eigenvector of \( \Q^{(k+1)} \) that is orthogonal to all the \( \phi_j^i \):s, i.e.\ \( \phi \) satisfies~\eqref{ee} and
\begin{equation} \label{ef}
-\Q^{(k+1)}\phi(x) = \lambda\phi(x)
\end{equation}
for all \( x=(v_1,\ldots,v_{k+1}) \in (V)_{k+1} \).
Spelling out the left hand side gives
\begin{equation} \label{eg}
\begin{split}
\left((k+1)(n-k-1) + {k+1 \choose 2}\right) \phi(x) - \sum_{\tau}\phi(x\tau) \\[-1ex] -
\sum_{i=1}^{k+1}\sum_{v\in V \colon v \not \in x}\phi(v_1,\ldots,v_{i-1},v,v_{i+1},\ldots,v_{k+1}),
\end{split}
\end{equation}
where \( \tau \) ranges over all \( {k+1 \choose 2} \) transpositions of two labelled balls and where we in the second term identified the \( k+1 \)-tuple \( x \) with its set of coordinates.
Using~\eqref{ee}, all the inner sums in the double sum in~\eqref{eg} simplifies to \( -\phi(x) \) and hence~\eqref{eg} simplifies to
\begin{equation} \label{eh}
\left((k+1)(n-k)+{k+1 \choose 2}\right) \phi(x) - \sum_{\tau}\phi(x\tau).
\end{equation}
Let \( H_x \) be the set of \( (k+1)! \) elements \( y \in (V)_{k+1} \) that one can get from \( x \) by permuting the labelled balls among themselves, but keeping the set of positions occupied by a labelled ball fixed.
Then considering~\eqref{eh} for \( y \in H_x \) and inserting in~\eqref{ef} becomes a ``local'' system of equations
 \begin{equation} \label{ei}
\left((k+1)(n-k)+{k+1 \choose 2}\right) \phi(y) - \sum_{\tau}\phi(y\tau) = \lambda \phi(y),
\end{equation}
\( y \in H_x \). This local system simply states that
\begin{equation} \label{ej}
\left(\left( (k+1)(n-k) + {k+1 \choose 2} \right) \I_{k+1} - \A_{k+1}\right) \phi|_{H_x} = \lambda\phi|_{H_x}
\end{equation}
To solve this, \( \lambda \) must be an eigenvalue of \( \left( (k+1)(n-k) + {k+1 \choose 2} \right) \I_{k+1} - \A_{k+1} \).
Since \( \spec(-\Q^{(k)}) \subseteq \cale^m_k \) by the induction hypothesis, this proves precisely that
\( \spec(-\Q^{(k+1)}) \subseteq \cale^m_{k+1} \) as desired.

To prove the claim of symmetry of \( \spec(\A_n) \) for any \( n \), note that \( \A \coloneqq \A_n \) is the adjacency matrix of a bipartite graph and can hence, by sorting the vertices of the graph appropriately,
be written in block form as
\[\A = \left[ \begin{array}{cc} 0 & \A_1 \\ \A_2 & 0 \end{array} \right]\]
However if \( \phi = [\phi_1 \,\, \phi_2]^T \) is an eigenvector of \( \A \) for the eigenvalue \( \lambda \), then \( [-\phi_1 \,\, \phi_2]^T \) is an eigenvector for \( -\lambda \),
proving the symmetry of the spectrum of \( \A \).

It remains to prove the multiplicity statements, that is, we need to prove that if \( 1 \leq j \leq \min{(\ell, n-\ell)} \), then the total multiplicity of the eigenvalues \( \lambda \) of \( -\Q^{(\ell)}\) that is such that
\[
j(n-j+1) \leq \lambda < (j+1)(n-j)
\]
is at most \( (n)_j \binom{\ell}{j} \).
To this end, note first that the total multiplicity of the eigenvectors of \( \Q^{(j)} \) can be at most \( |(V)_j| = (n)_j \).
Secondly, note that for any \( j \), the largest eigenvalue of \( \mathbb{A}_{j} \) is \(\binom{j}{2} \). Using~\eqref{ej}, it follows that any \emph{new} eigenvalue \( \lambda \)  we get on level \( j \), i.e. an eigenvalue that does not correspond to a lifted eigenvector, satisfies
\begin{equation}\label{ma}
\lambda \geq \left(  j(n-j+1) + \binom{j}{2} \right) - \binom{j}{2} = j(n-j+1).
\end{equation}
Consequently, it now follows that any eigenvalue \( \lambda' \)  of \( -\Q^{(\ell)}\) that is such that
\[
\lambda' < (j+1)(n-j)
\]
must correspond to a lifted eigenvector from either level at most \( j \) or level at least \(n-j+1\).
If the second holds, we must have that \( n-j+1 \leq \ell \), or equivalently, that \( n - \ell <    j  \), which contradicts that \( j \leq \min{(\ell, n-\ell) }\), so the first of these must hold, that is \( \lambda' \) must correspond to  a lifted eigenvector of \( -\Q^{(j)} \). The number of eigenvectors of \( -\Q^{(j)} \) is exactly \( (n)_j\), and these can be lifted in at most \( \binom{\ell}{j} \) ways, why the desired conclusion follows.

%It remains to prove the multiplicity statements. Since \( |(V)_k| = (n)_k \), the total multiplicity of eigenvalues \( \lambda \) with \( k(n-k+1) \leq \lambda < (k+1)(n-k) \) of \( \Q^{(n,k)} \) can obviously not exceed \( (n)_k \). However as long as \( k < n/2 \), \( k(n-\ell+1) \) is increasing in \( \ell \), so by what we found above, adding extra labelled particles does not result in any new eigenvalues below \( (k+1)(n-k) \) and the eigenvectors of \( -\Q^{(n,\ell)} \), \( k<\ell \leq n/2 \), corresponding to eigenvalues of \( -\Q^{(n,k)} \) can only depend on at most \( k \) labelled balls. Since there are less than \( (n)_k \) linearly independent such eigenvectors for \( k \) given labelled balls, having \( \ell \) balls from which one can choose \( {\ell \choose k} \) subsets of size \( k \), it follows that the total multiplicity of eigenvalues with \( k(n-k+1) \leq \lambda < (k+1)(n-k) \) as eigenvalues of \( -\Q^{(n,\ell)} \) is at most \( {\ell \choose k}(n)_k \).

For the final claim that the multiplicity of the eigenvalue \( \lambda=n \) of \( -\Q^{(\ell)} \) is \( \ell(n-1) \) follows from a simplified version of this argument: for one given ball we know from the UEP that
the multiplicity is \( n-1 \) and there are thus \( n-1 \) orthogonal eigenvectors \( \phi_1,\ldots,\phi_{n-1} \).
Then, when we have \( \ell \) labelled balls to choose from, we define for each ball \( i \), \( \phi^i_j(v_1,\ldots,v_\ell) = \phi_j(v_i) \).
As before, \( \phi^i_j \) and \( \phi^i_{j'} \) are orthogonal for \( j \neq j' \).
Since \( \psi_1 \), \ldots, \( \psi_{n-1} \) are all orthogonal to \( (1,1,\ldots, 1) \), we obtain that \( \sum_{v=1}^n\phi_j(v)=0 \) for all \( j \) and  it also follows that \( \phi^i_j \) and \( \phi^l_{j'} \) are orthogonal for \( i \neq \ell \) for any \( (j,j') \).
It follows that \( \{\phi^i_j\}_{i=1,\ldots,\ell,\,j=1,\ldots,n-1} \) is an orthogonal family of eigenvectors. To prove the claim, we now only need to argue that there can be no eigenvectors that are orthogonal to these vectors with the same eigenvalue. However, from~\eqref{ma} it follows that any such eigenvector must be lifted from level at most 1. As these have already been considered, the desired conclusion follows.
\end{proof}

\section{Proofs of \( L^2 \)-mixing times}\label{sec:4}

Consider an irreducible continuous time Markov chain \( \{X_t\}_{t \geq 0} \) on a finite state space \( S \) with a symmetric generator \( \Q \).
Since \( \Q \) is symmetric, the stationary distribution \( \pi \) is uniform. Let \( N\coloneqq|S| \). In this section we will let \( \langle \cdot , \cdot \rangle \) be the usual inner product on \( \mathbb{R}^N \):
\[ \langle f,g \rangle = \sum_{s \in S}f(s)g(s).\]
Note that as \( \pi  \) is uniform, this inner product differs from the inner product used earlier in this paper only by a scaling.
Let \( 0=\mu_1<\mu_2 \leq \ldots \leq \mu_N \) be the eigenvalues of \( \Q \) and let \( \phi_1,\ldots,\phi_N \) be an orthonormal family of corresponding eigenvectors.
Let \( x \in S \) be the starting state of \( \{ X_t\}_{t \geq 0}\) and write the function \( e_{x}(s)=\one_{x}(s) \) in the eigenvector basis as
\[e_{x} = \sum_{i=1}^N c_i(x) \phi_i,\]
where \( c_i(x)=\langle e_x,\phi_i \rangle \).
Note that as with this scaling, \( \phi_1 \equiv 1 / \sqrt{|S|}\) and
\[
c_1(x) = \langle e_x, \phi_1 \rangle = \sum_{s \in S} e_x(s) \phi_1(x) = \phi_1(x) = 1/\sqrt{|S|}
\]
the term \( c_1(x) \phi_1 = 1/\sqrt{|S|} \cdot 1/\sqrt{|S|} = 1/|S| = \pi(x) \).
Hence by standard arguments
\[\Pro_x(X_t \in \cdot)-\pi = \sum_{i=2}^N c_i(x)e^{-\mu_i t}\phi_i\]
and consequently
\[\n \Pro_x(X_t \in \cdot) - \pi \n_2^2 = N\sum_{i=2}^N c_i(x)^2e^{-2\mu_i t}.\]
Let us now write \( 0=\lambda_1 < \lambda_2 < \ldots < \lambda_r=\mu_N \) for the {\em distinct} eigenvalues of \( \Q \) (so \( r \leq N \)).
For each \( j=1,\ldots,r \), let \( C_j(x)^2 = \sum_{i\colon \mu_i = \lambda_j}c_i(x)^2 \).
Then we can rewrite as
\begin{equation} \label{ek}
\n\Pro_x(X_t \in \cdot) - \pi \n_2^2= N\sum_{j=2}^r C_j(x)^2e^{-2\lambda_j t}.
\end{equation}
Now assume that our Markov chain is such that \( \n \Pro_x(X_t \in \cdot)-\pi\n_2 \) is independent of \( x \), such as is the case for the UEP and the LEP.
It then follows that \( C_j(x) \) is independent of \( x \).
Let \( m_j \) be the multiplicity of the eigenvalue \( \lambda_j \).

\bla \label{le}
If \( \{X_t\} \) is such that \( \n \Pro_x(X_t \in \cdot)-\pi\n_2 \) is independent of \( x \), then for every \( j \) and every \( x \),
\[C_j(x)^2=\frac{m_j}{N}.\]
\ela

\bpf

For every \( x \in S \), let \( e_x \in \RR^N \) be the corresponding unit vector. Fix \( j \in \{ 1,2, \ldots, r \} \) and note that the eigenvectors for the eigenvalue \( \lambda_j \) span a subspace \( U_j\) of dimension \( m_j \). By symmetry, we know that the projection of each \( e_x \) onto \( U_j \) has the same length.
%We have the unit vectors \( e_j \in \RR^N \), \( j \in S \), a subspace \( U \) of dimension \( m \) and we know that the projection of each \( e_j \) onto \( U \) has the same length.

Let \( u_i = (u_{i1},\ldots,u_{iN}) \), \( i=1,\ldots,m_j \) be an orthonormal basis for \( U_j \).
The projection of \( e_x \) onto \( u_i \) is \( u_{ix} \), so the square length, \( C_j(x)^2 \), of the projection of \( e_x \) onto \( U_j \) is \( \sum_{i=1}^{m_j} u_{ix}^2 \).
Summing over \( x \) gives
\[\sum_{x=1}^N C_j(x)^2 = \sum_{i=1}^{m_j} \sum_{j=1}^N u_{ij}^2 = m_j\]
since the \( u_i \):s are unit vectors.
Since the \( C_j(x)^2 \):s are equal, they must all equal \( m_j/N \).
This proves the lemma.
\epf

Applying Lemma~\ref{le} to~\eqref{ek}, it follows that in situations where the \( L^2 \)-norm does not depend on the starting state,
\begin{equation} \label{el}
\n\Pro(X_t \in \cdot)-\pi\n_2^2 = \sum_{j=2}^r m_j e^{-2\lambda_j t}.
\end{equation}
This together with Theorem~\ref{ta}, recalling that \( m_j = \binom{n}{j-1} - \binom{n}{j-2}\) and that \( \alpha=2/n^2 \) so that the eigenvalues are \( \lambda_j = 2(j-1)(n-j+2)/n^2 \),
proves the exact formula for the \( L^2 \)-distance of Theorem~\ref{tc}.
It remains to check the estimates.
For the lower bound, it suffices to recall that  \( t=(1/4)n\log(n-1)+cn \) and that the first term is
\[(n-1)e^{-4t/n} = e^{-4c}.\]
Taking the square root gives the result.

For the upper bound, take \( c \geq 0 \) and observe that
\begin{align*}
\n\Pro(X_t \in \cdot)-\pi\n_2^2 &= \sum_{j=2}^r m_j e^{-2\lambda_j t} \\
&< \sum_{j=2}^\ell {n \choose j-1}e^{-(j-1)(n-j+2)(\log n + 4c)/n} \\
&= \sum_{j=1}^{\ell-1} {n \choose j}e^{-j(n-j+1)(\log n + 4c)/n} \\
&< \sum_{j=1}^{n/2} \frac{n^j}{j!}e^{-j\log n}e^{j(j-1)\log n/n}e^{-4j(n-j+1)c/n} \\
&< e^{-4c}\sum_{j=1}^{n/2}\frac{n^{j(j-1)/n}}{j!}.
\end{align*}
Let \( n \geq 1000 \) and \( 10 \leq j \leq n/2 \).
Taking logarithms and using the estimate \( \log j! \geq j \log j - j\), it is easy to see that for such \( n \) and \( j \), \( n^{j(j-1)/n}/j! < e^{-j} \).
Hence
\[\sum_{j=10}^{n/2}\frac{n^{j(j-1)/n}}{j!} < \sum_{j=10}^\iy e^{-j} < e^{-9}.\]
Also, for \( j<10 \) and \( n \geq 1000 \), we have \( j<n^{1/3} \), so for \( n \geq 1000 \),
\[\sum_{j=1}^{9}\frac{n^{j(j-1)/n}}{j!} < n^{1/n^{1/3}}\sum_{j=1}^\iy\frac{1}{j!} < 2(e-1).\]
Summing up gives
\[\n\Pro(X_t \in \cdot)-\pi\n_2^2 < (e^{-9}+2(e-1))e^{-4c} < 4e^{-4c}.\]
Now take square roots again to finish the proof of Theorem~\ref{tc}.

\bigskip

Let us now move to the LEP. For the lower bound, we again simply consider the first term of the the right hand side of~\eqref{el}.
By Theorem~\ref{tb}, the multiplicity \( m_2 \) is \( \ell(n-1) \) and \( \lambda_2 \) is still \( 2/n \), so the first term now becomes
\[\ell(n-1)e^{-4t/n} = e^{-4c},\]
using that \( t \) is now \( (1/4)n\log(\ell(n-1))+cn \).
Taking square roots gives the desired lower bound.

For the upper bound, take \( c \geq 0 \). If \( \ell \leq (1-\varepsilon)n/2 \) for some \( \varepsilon > 0 \), using the multiplicity bounds of Theorem~\ref{tb}, that all eigenvalues are positive and  that \( t =  (1/4)n\log(\ell(n-1))+cn \), we find that
\begin{align*}
\n\Pro(X_t \in \cdot)-\pi\n_2^2 &\leq \sum_{j=1}^{\ell} {\ell \choose j}(n)_j e^{-j(n-j+1)(\log(\ell n)+4c)/n}  \\
&<  \sum_{j = 1}^\ell  \frac{(\ell n)^j}{j!}e^{-j\log(\ell n)}e^{j(j-1)/n\cdot \log(\ell n) }e^{-4j(n-j+1)c/n}   \\
&< e^{-4c}\sum_{j=1}^\ell \frac{(\ell n)^{j(j-1)/n}}{j!}  \\
&< e^{-4c}\sum_{j=1}^{(1-\varepsilon)n/2} \frac{n^{2j(j-1)/n}}{j!} .
\end{align*}

Taking logarithms and using Stirling's formula, it follows that whenever \( {n \geq \max\left(8000,\exp\left( (\log 2 + 1 - \log(1 - \varepsilon))/\varepsilon\right)\right)} \) and \( 20 \leq j \), we have \( n^{2j(j-1)/n}/j! < e^{-j} \).
Hence
\[\sum_{j=20}^{\ell} \frac{n^{2j(j-1)/n}}{j!} < \sum_{j=20}^\iy e^{-j} < e^{-19}.\]
For \( j < 20 \) and \( n \geq 8000 \), we have \( j \leq n^{1/3} \), so for \(  {n \geq \exp\left( (\log 2 + 1 - \log(1 - \varepsilon))/\varepsilon\right)} \),
\[\sum_{j=1}^{19}\frac{n^{2j(j-1)/n}}{j!} < n^{2/n^{1/3}}\sum_{j=1}^\iy\frac{1}{j!} < 2(e-1).\]
Summing up gives
\[\n\Pro(X_t \in \cdot)-\pi\n_2^2 < (e^{-19}+2(e-1))e^{-4c} \leq 4e^{-4c}.\]
This establishes the upper bound of Theorem~\ref{td}.

\end{document}